\documentclass[11pt]{article}

\usepackage{amssymb,latexsym,amsmath,amsthm, mathrsfs}
\usepackage{graphicx}
\usepackage{setspace}
\usepackage[all]{xy}
\newcommand{\susp}{\Sigma\mkern-10.1mu/}

\newtheorem{theorem}{Theorem}[section]
\newtheorem{lemma}[theorem]{Lemma}
\newtheorem{proposition}[theorem]{Proposition}
\newtheorem{corollary}[theorem]{Corollary}

\newenvironment{definition}[1][Definition]{\begin{trivlist}
\item[\hskip \labelsep {\bfseries #1}]}{\end{trivlist}}
\newenvironment{example}[1][Example]{\begin{trivlist}
\item[\hskip \labelsep {\bfseries #1}]}{\end{trivlist}}
\newenvironment{remark}[1][Remark]{\begin{trivlist}
\item[\hskip \labelsep {\bfseries #1}]}{\end{trivlist}}
\newenvironment{acknowledgement}[1][Acknowledgements]
{\begin{trivlist} \item[\hskip \labelsep {\bfseries #1}]}
{\end{trivlist}}

\DeclareMathSymbol{\N}{\mathbin}{AMSb}{"4E}
\DeclareMathSymbol{\Z}{\mathbin}{AMSb}{"5A}
\DeclareMathSymbol{\R}{\mathbin}{AMSb}{"52}
\DeclareMathSymbol{\Q}{\mathbin}{AMSb}{"51}
\DeclareMathSymbol{\I}{\mathbin}{AMSb}{"49}
\DeclareMathSymbol{\C}{\mathbin}{AMSb}{"43}

\def\P{{\mathbb P}}
\def\RP{{\mathbb R \mathbb P}}
\def\CP{{\mathbb C \mathbb P}}

\title{Complexification of real cycles and Lawson suspension theorem}
\author{Jyh-Haur Teh\\Department of Mathematics
\\National Tsing Hua University of Taiwan}
\begin{document}
\date{}
\maketitle

\footnotetext[1]{Tel: +886 3 5715131-33071\\
E-mail address: jyhhaur@math.nthu.edu.tw}

\begin{abstract}
We embed the space of totally real $r$-cycles of a totally real
projective variety into the space of complex $r$-cycles by
complexification. We provide a proof of the holomorphic taffy
argument in the proof of Lawson suspension theorem by using Chow
forms and this proof gives us an analogous result for totally real
cycle spaces. We use Sturm theorem to derive a criterion for a real
polynomial of degree $d$ to have $d$ distinct real roots and use it
to prove the openness of some subsets of real divisors. This enables
us to prove that the suspension map induces a weak homotopy
equivalence between two enlarged spaces of totally real cycle
spaces.
\end{abstract}

\section{Introduction}
Classically, real algebraic geometry is the study of complex
varieties endowed with a real structure. For example, with complex
conjugation, the complex projective space $\CP^n$ is a real
variety. But there is another candidate for the name real
algebraic geometry which is the study of real zero loci of real
polynomials. The development in this branch of geometry shows that
it has many interesting problems and people start indicating it
the real algebraic geometry, for example in \cite{AK}, \cite{BCR}.
Since we deal with both of these geometries in this paper, to
distinguish them, we will keep the name real algebraic geometry
for the classical involutional geometry and totally real algebraic
geometry for the other one.

This new branch sits between differential topology and algebraic
geometry. By the apparatus of Morse theory, Milnor (see \cite{M})
and Thom (see \cite{Thom}) gave an estimation of the total Betti
number of a totally real algebraic variety. This estimation has been
applied for instance in obtaining lower bounds on the complexity of
algorithms (see \cite{Ben}). A fascinating fundamental result in
totally real algebraic geometry is the Nash-Tognoli theorem which
says that every compact smooth manifold without boundary is
diffeomorphic to some nonsingular algebraic varieties. Since complex
projective manifolds are more restricted, there is no analogous
result in complex algebraic geometry.

In general, totally real algebraic varieties are more irregular than
real algebraic varieties. One of the most striking differences is
that the real projective space $\RP^n$ is an affine variety. From
here it is not difficult to see that it is difficult to study
totally real algebraic varieties from the classical methods like
intersection or projection. Nevertheless, many properties of a
totally real projective variety behave well in $\Z_2$-world. We
believe that a result in complex world should have its analogue in
the real world but with $\Z_2$-nature. For example the classical
Fundamental Theorem of Algebra is certainly not true in the real
world, but we have a $\Z_2$-version of it. Another example is Lawson
homology theory and morphic cohomology theory for projective
varieties founded by Lawson and Friedlander. The parallel theories
for real projective varieties are called reduced real Lawson
homology and reduced real morphic cohomology developed by the author
(see \cite{T1}). The cornerstone of all these theories is Lawson
suspension theorem (see \cite{F}, \cite{L}) which says that the
suspension map induces a homotopy equivalence between $Z_r(X)$ and
$Z_{r+1}(\susp X)$ for a projective variety $X$ where $Z_r(X)$ is
the $r$-cycle group of $X$. As an application of this theorem, we
have
$$Z_r(\P^n)\cong K(\Z, 0)\times K(\Z, 2)\times \cdots \times K(\Z, 2(n-r))$$
where $K(\Z, i)$ is the Eilenberg-Mac lane space. A version of the
Lawson suspension theorem in reduced real Lawson homology is the
reduced real Lawson suspension theorem which says that $R_r(X)$ is
homotopy equivalent to $R_{r+1}(\susp X)$ for a real projective
variety $X$ where $R_r(X)$ is the reduced real cycle group of $X$.
As an application, we have
$$R_r(\P^n)\cong K(\Z_2, 0)\times K(\Z_2, 1)\times \cdots \times K(\Z_2,
n-r).$$ We note that the real points of the complex projective
space $\P^n$ is $\RP^n$ and the nonzero homology groups
$H_i(\RP^n;\Z_2)$ with $\Z_2$-coefficients of $\RP^n$ are $\Z_2$
for $i=1,..., n$.

The reduced real Lawson homology and the reduced real morphic
cohomology enable us to use real projective varieties to study
totally real projective varieties. One reason is that there are
natural transformations from the reduced real Lawson homology and
the reduced real morphic cohomology of a real projective variety to
the singular homology groups and the singular cohomology groups with
$\Z_2$-coefficients respectively of its real points. For instance, a
generalization of the Harnack-Thom theorem is proved in \cite{T2}
which relates the rank of the Lawson homology groups with
$\Z_2$-coefficients of a real projective variety $X$ with the rank
of its reduced real Lawson homology groups. For zero-cycles, we
recover the classical Harnack-Thom theorem which gives a bound of
the total Betti number in $\Z_2$-coefficients of the real points by
the total Betti number in $\Z_2$-coefficients of a real projective
variety.

In this paper, we study the space formed by totally real $r$-cycles
on a totally real projective variety. A totally real $r$-cycle in
our sense is a formal sum of some $r$ dimensional totally real
irreducible subvarieties on a totally real variety. We embed the
space of totally real $r$-cycles of a totally real variety $X_{\R}$
into the space of complex $r$-cycles of its complexification
$X_{\C}$. Therefore it inherits the topology from the complex
$r$-cycle group. Our goal is to prove a result analogous to the
Lawson suspension theorem. The proof of Lawson suspension theorem is
divided into two parts which Lawson calls \emph{holomorphic taffy}
and \emph{magic fans}. Since the space of totally real $r$-cycles is
not closed in the \emph{holomorphic taffy} process, we have to
enlarge the embedded space of totally real $r$-cycles $RZ_r(X)$ on
$X$ to a space $SZ_r(X)$ and enlarge the embedded space of totally
$(r+1)$-cycles $RZ_{r+1}(\susp X)$ on $\susp X$ to a space
$\overline{RZ}_{r+1}(\susp X)$. We follow his division to show that
the suspension map induces a weak homotopy equivalence between
$SZ_r(X)$ and $\overline{RZ}_{r+1}(\susp X)$. We use this suspension
theorem to compute the homotopy type of $\overline{RZ}_1(\RP^n)$.
The algebraic proof of Lawson suspension theorem given by
Friedlander in \cite{F} works for any algebraically closed fields.
His proof is schematic and needs more machinery. Our proof of the
part \emph{holomorphic taffy} deals with Chow forms directly and is
easier to apply to the proof of \emph{holomorphic taffy} for totally
real cycles.

\section{Complexification}
In this paper, a totally real algebraic variety is the zero loci
in $\R^n$ of some real polynomials and a totally real projective
variety is the zero loci in $\RP^n$ of some real homogeneous
polynomials. Throughout this paper, $\R^n$ and $\RP^n$ will be
considered as subsets of $\C^n$ and $\CP^n$ respectively under the
canonical embeddings. For the basic properties of totally real
algebraic varieties we refer to \cite{AK}, \cite{BCR}.

\begin{definition}
Suppose that $X \subseteq \RP^n$ is a subset. The real cone
$RC(X)$ of $X$ is $\pi^{-1}(X)\cup \{ 0\}$ where $\pi$ is the
quotient map $\pi: \R^{n+1}-{0} \rightarrow \RP^n$. We denote the
ideal in $\R[x_0,..., x_n]$ generated by those real homogeneous
polynomials vanishing on $RC(X)$ by $I^{\R}(X)$. Now consider
$RC(X)$ as a subset in $\C^{n+1}$ and we denote $I^{\C}(X)$ to be
the ideal in $\C[z_0,...,z_n]$ generated by those homogenous
polynomials vanishing on $RC(X)$. For a homogeneous ideal $J
\subseteq \C[z_0,...,z_n]$, we denote the complex projective
variety in $\P^n$ defined by $J$ by $Z(J)$. If $X$ is a totally
real projective variety, its complexification $X_{\C}$ is defined
to be $ZI^{\C}(X)$.
\end{definition}

\begin{definition}
For $f \in \R[x_0,...,x_n]$, the complexification of $f$ is
defined to be the complex polynomial $f_{\C} \in \C[z_0,...,z_n]$
where $f_{\C}$ is the extension of $f$ to $\C^{n+1}$. In other
words, it is the image of $f$ under the natural embedding
$\R[x_0,..., x_n] \hookrightarrow \C[z_0,..., z_n]$. For an ideal
$J \subset \R [x_0,...,x_n]$, the complexification $J_{\C}$ of $J$
is the ideal generated by $\{f_{\C}| f \in J\}$ in
$\C[z_0,...,z_n]$.
\end{definition}

\begin{proposition}
For a totally real projective variety $X\subset \RP^n$,
$I^{\C}(X)=(I^{\R}(X))_{\C}.$ \label{complexifying ideal}
\end{proposition}

\begin{proof}
By definition, if $f \in I^{\C}(X)$, then $f(RC(X))=0$. Let $g$ be
the restriction of $f$ to $\R^{n+1}$, and let $Re(g), Im(g)$ be
two real polynomials such that $g=Re(g)+iIm(g)$. Since
$g(RC(X))=f(RC(X))=0$, $Re(g)(RC(X))=Im(g)(RC(X))=0$. So $Re(g),
Im(g)$ are in $I^{\R}(X)$ which implies that $Re(g)_{\C}$,
$Im(g)_{\C}$ are in $(I^{\R}(X))_{\C}$. So
$f=Re(g)_{\C}+iIm(g)_{\C} \in (I^{\R}(X))_{\C}$. Another direction
follows from the definition.
\end{proof}

\begin{definition}
Suppose that $V \subset \R^n$ is a totally real algebraic variety
and $I(V)$ is the ideal of polynomials in $\R^n$ vanishing on $V$.
The dimension $dimV$ of $V$ is defined to be the maximal length of
chain of prime ideals in $\R[x_1,..., x_n]/I(V)$. The dimension of
a totally real projective variety is the dimension of its real
cone minus 1. We say that a point $x \in V$ is nonsingular if and
only if we can find a generator set $\{p_1,..., p_k\}$ of $I(V)$
such that
$$rk[\frac{\partial p_i}{\partial x_j}(x)]= n-dimV.$$ A point
of a totally real projective variety is called nonsingular if it
is the image of a nonsingular point in the real cone of that
totally real projective variety.
\end{definition}

The following is from \cite{BCR}, Proposition 3.3.8.
\begin{proposition}
If $V \subset \R^n$ is a totally real algebraic variety of
dimension $d$ and $x \in V$ is nonsingular. Then there exist
$p_1,..., p_{n-d} \in I(V)$ and a Zariski open neighborhood $U$ of
$x$ in $\R^n$ such that

i. $Z^{\R}(p_1,..., p_{n-d}) \cap U=V\cap U$ where
$Z^{\R}(p_1,..., p_{n-d})$ is the totally real algebraic variety
defined by $p_1,..., p_{n-d}$.

ii. for every $y \in U$, rank $[\frac{\partial p_i}{\partial
x_j}(y)]=n-d.$
\end{proposition}

The following is the key proposition which enables us to embed
spaces of totally real cycles into complex cycle spaces. To
simplify the notation, we use $ZI(X)$ to designate $ZI^{\C}(X)$.
We note that for a totally real projective variety $X\subset
\RP^n$, $ZI(X)$ is the smallest complex projective variety in
$\P^n$ which contains $X$.

\begin{proposition}
Suppose that $X \subset \RP^n$ is a totally real projective
variety.
\begin{enumerate}
\item If $Y \subset X$ is a totally real projective subvariety,
then $ZI(Y) \cap X=Y$.

\item If $Y$ is an irreducible p-dimensional totally real
projective subvariety of $X$, then $ZI(Y)$ is an irreducible
complex p-dimensional projective subvariety of $ZI(X)$.

\item If $V$ is a complex p-dimensional irreducible projective
subvariety of $ZI(X)$ and dim$_{\R}(V \cap X)=p$, then $V \cap X$
is an irreducible totally real projective subvariety of $X$ and
$ZI(V \cap X)=V$.
\end{enumerate}
\label{complexification}
\end{proposition}

\begin{proof}
\begin{enumerate}

\item Let $q\in ZI(Y) \cap X$. For any $f \in I^{\R}(Y)$, $f_{\C}
\in I^{\C}(Y)$, thus $f_{\C}(q)=0$. But $q$ is a real point, so
$f(q)=0$. Therefore $ZI(Y) \cap X$ is contained in $Y$. Another
direction is trivial.

\item Write $ZI(Y)=V_1 \cup \cdots \cup V_t$ as a union of its
irreducible components. Let $Y_i=Y \cap V_i$ for each $i$. Then we
have $\cup Y_i=Y$. But $Y$ is irreducible, so some $Y_i$ equals to
$Y$. This means that $V_i$ contains $Y$. But $ZI(Y)$ is the
smallest complex projective variety contains $Y$, thus
$V_i=ZI(Y)$. So $ZI(Y)$ is irreducible. Now let us consider the
dimension of $ZI(Y)$. Let $NS$ and $Sing$ be the nonsingular part
and singular part of a variety respectively. It is known that the
singular part of a projective variety is again a projective
variety. We claim that the nonsingular part of $Y$ and the
nonsingular part of the complexification of $Y$ must intersect
somewhere. If not, then $NS(Y) \subset Sing(ZI(Y))$. Take
$W=ZI(NS(Y)) \subset Sing(ZI(Y))$. Since $Y=(W \cap Y) \cup
Sing(Y)$ and $Y$ is irreducible, we have $Y=W \cap Y$. But $ZI(Y)$
is the smallest complex projective variety which contains $Y$, we
must have $ZI(Y) \subset W \subset Sing(ZI(Y))$. A contradiction.
Now take a point $q$ which is a nonsingular point of $Y$ and
$ZI(Y)$ and assume that $I^{\R}(Y)$ is generated by $f_1,...,
f_k$. By Proposition \ref{complexifying ideal}, we know that
$I^{\C}(Y)$ is the complexification of $I^{\R}(Y)$ and therefore
$I^{\C}(Y)$ is generated by $f_{1, \C},..., f_{k, \C}$. Since $q$
is nonsingular, thus
$$n-dim_{\R}(Y)=rk[\frac{\partial f_i}{\partial x_j}](q)=
rk[\frac{\partial f_{i \C }}{\partial z_j}](q)=n-dim_{\C}ZI(Y).$$
We have $dim_{\R}(Y)=dim_{\C}(ZI(Y))$.

\item Suppose that $Y\subset V\cap X$ is an irreducible totally
real projective subvariety of dimension $p$. Then $ZI(Y) \subset
V$. By (2), dim$_{\C}ZI(Y)=p$ and $V$ is irreducible, so
$ZI(Y)=V$. By (1), $Y=ZI(Y) \cap X=V \cap X$. So $V \cap X$ is
irreducible and of dimension $p$. By (2), $ZI(V \cap X)$ is of
dimension $p$ and $ZI(V \cap X)\subset V$. The irreducibility and
the dimension of $V$ imply that $ZI(V \cap X)=V$.
\end{enumerate}
\end{proof}

\begin{definition}
Suppose that $X \subset \RP^n$ is a totally real projective
variety. A totally real $r$-cycle on $X$ is a linear combination
of irreducible totally real subvarieties of dimension $r$ with
integral coefficients. We say that a totally real cycle is
effective if all of its coefficients are positive. The
complexification of a totally real cycle $c=\sum_{i=0}^nn_iV_i$ is
defined to be the complex cycle
$\widetilde{c}=\sum_{i=0}^nn_iZI(V_i)$ on $X_{\C}$ and the degree
of a totally real cycle is the degree of its complexification.
\end{definition}

\begin{remark}
The geometric picture of an irreducible totally real subvariety is
very different from the complex one. For instance, an irreducible
totally real algebraic variety may not be connected which
contrasts to the case of complex irreducible algebraic variety. A
concrete example is the irreducible cubic curve in $\R^2$ given by
the equation $x^2+y^2-x^3=0$ which has two connected components
and furthermore, two connected components have different
dimensions.
\end{remark}

\begin{definition}
Suppose that $X \subset \RP^n$ is a totally real projective
variety. Let $R\mathscr{C}_{p,d}(X)$ be the collection of all
effective totally real $p$-cycles of degree $d$ on $X$.
\end{definition}

We denote the Chow variety of effective $p$-cycles of degree $d$
on $X_{\C}$ by $\mathscr{C}_{p, d}(X_{\C})$ and denote the Chow
monoid of $p$-cycles by $\mathscr{C}_p(X_{\C})=\underset{d\geq
0}{\coprod}\mathscr{C}_{p, d}(X_{\C})$. Let
$$K_{p, d}(X_{\C})=\underset{d_1+d_2\leq
d}{\coprod}\mathscr{C}_{p, d_1}(X_{\C})\times \mathscr{C}_{p,
d_2}(X_{\C})/\sim$$ where the equivalence relation $\sim$ is
defined by $(a, b) \sim (c, d)$ if and only if $a+d=b+c$. We give
each $K_{p, d}(X_{\C})$ the quotient topology and since each Chow
variety is compact, $K_{p, d}(X_{\C})$ is compact for all $d$. The
group of $p$-cycles $Z_p(X_{\C})$ is the Grothendieck group (naive
group completion) of $\mathscr{C}_p(X_{\C})$ and from the
filtration
$$K_{p, 0}(X_{\C})\subset K_{p, 1}(X_{\C}) \subset
\cdots = Z_p(X_{\C}),$$  we give $Z_p(X_{\C})$ the weak topology
defined by this filtration which means that a set $W \subset
Z_p(X_{\C})$ is closed if and only if $W\cap K_{p, d}(X_{\C})$ is
closed for all d.

From Proposition \ref{complexification}, we see that
$R\mathscr{C}_{p, d}(X)$ is mapped injectively into
$\mathscr{C}_{p, d}(X_{\C})$ by taking the complexification of
totally real cycles. So from now on we identify $R\mathscr{C}_{p,
d}(X)$ with its image in $\mathscr{C}_{p, d}(X_{\C})$. Let
$$R\mathscr{C}_p(X)=\underset{d\geq 0}{\coprod}R\mathscr{C}_{p,
d}(X)$$ be the monoid formed by totally real $p$-cycles and denote
its naive group completion by $RZ_p(X)$. Endow $RZ_p(X)$ with the
subspace topology of $Z_p(X)$, then $RZ_p(X)$ is a topological
subgroup of $Z_p(X_{\C})$. We denote the closure of $RZ_p(X)$ in
$Z_p(X_{\C})$ by $\overline{RZ}_p(X)$.

\section{A Proof of the Holomorphic Taffy by Chow Forms}
Suppose that $V^r \subset \P^n$ is a $r$-dimensional irreducible
subvariety. For a hyperplane $\xi$ in $\P^n$, it corresponds to a
point $\underline{\xi}$ in $\P^n$ by taking its coefficients.
Obviously this correspondence is a bijection between the set of
hyperplanes in $\P^n$ and $\P^n$. Let $C$ be the set consisting of
$(\underline{a}, \underline{\xi}^0,...,\underline{\xi^r}) \in V
\times \P^n \times \cdots \times \P^n$ where $\xi^0,..., \xi^r$
are hyperplanes in $\P^n$ which meet $V$ at $a$. The set $C$ is an
irreducible projective variety whose ideal is generated by
$$I_n(V), \mbox{ } \sum_{j=0}^n\xi_j^ix_j, \forall i=0,...,r$$
where $\underline{a}=(a_0:...:a_n)$ and
$\underline{\xi}^i=(\xi_0^i:...:\xi_n^i)$. Project $C$ to
$(\P^n)^{r+1}$ we get a hypersurface. Thus there is an irreducible
multihomogenous form $F_V(\underline{u}^0,...,\underline{u}^r)$,
unique up to a constant, defines that hypersurface. The
multihomogeneous form $F_V$ is called the Chow form of $V$. From
here we can see that:

\begin{proposition}
An irreducible multihomogeneous polynomial
$F(\underline{u}^0,...,\underline{u}^r)$ of  $(\P^n)^{r+1}$ is the
Chow form of an irreducible subvariety $V \subset \P^n$ if and only
if it satisfies the following property: for any (r+1) hyperplanes
$\xi^0,...,\xi^r$, $F(\underline{\xi}^0,...,\underline{\xi}^r)=0$ if
and only if $\xi^0 \cap ...\cap\xi^r\cap V \neq \emptyset$. Or
equivalently, for $a\in \P^n$, $a \in V$ if and only if for any
$(r+1)$ hyperplanes $\xi^0,...,\xi^r$ which contain $a$,
$F(\underline{\xi}^0,...,\underline{\xi}^r)=0$.
\end{proposition}

Fix $t \in \C \backslash \{0\}$. Define $\varphi_t:\P^{m+n+1}
\rightarrow \P^{m+n+1}$ by
$$\varphi_t(y_0:...:y_m:x_0:...:x_n)=(ty_0:...:ty_m:x_0:...:x_n).$$
For $V$ an irreducible subvariety of $\P^n$, it is easy to see
that $\varphi_t(V)$ is also an irreducible subvariety of $\P^n$.

Let us define some notation. For $\underline{u}=(u_0:...:
u_{m+n+1})\in \P^{m+n+1}$, denote
$^t\underline{u}=(tu_0:...:tu_m:u_{m+1}:...:u_{m+n+1})$. And for a
hyperplane $\xi$ defined by $\sum^{m+n+1}_{i=0}\xi_ix_i=0$,
$^t\xi$ is the hyperplane defined by
$\sum^m_{i=0}t\xi_ix_i+\sum^{m+n+1}_{i=m+1}\xi_ix_i=0$.

\begin{lemma}
Suppose that $F_V$ is the Chow form of a $r$-dimensional
irreducible subvariety $V \subset \P^n$. Then
$^tF_V(\underline{u}^0,...,\underline{u}^r):=
F_V(^t\underline{u}^0,...,^t\underline{u}^r)$ is the Chow form of
$\varphi_t(V)$.
\end{lemma}

\begin{proof}
Suppose that $a\in \varphi_t(V)$ and $\xi^0,...,\xi^r$ are
hyperplanes in $\P^n$ which contain $a$. Then $^{\frac{1}{t}}a \in
V$ and $^t\xi^0,...,^t\xi^r$ contain $^\frac{1}{t}a$. Thus
$F_V(^t\underline{\xi}^0,..., ^t\underline{\xi}^r)=0$ which
implies $^tF_V(\underline{\xi}^0,..., \underline{\xi}^r)=0$. Now
we assume that for any $(r+1)$ hyperplanes which contain $a$,
$^tF_V$ evaluates at the point defined by these hyperplanes is
$0$. We need to show that $a$ is in $\varphi_t(V)$. Suppose that
$\xi^0,...,\xi^r$ are any $(r+1)$ hyperplanes which contain
$^{\frac{1}{t}}a$. Then $^\frac{1}{t}\xi^0,....,^\frac{1}{t}\xi^r$
are hyperplanes which contain $a$. Thus
$^tF_V(^\frac{1}{t}\underline{\xi}^0,...,^\frac{1}{t}\underline{\xi}^r)=
F_V(\underline{\xi}^0,...,\underline{\xi}^r)=0$. So
$^{\frac{1}{t}}a \in V$ which implies $a \in \varphi_t(V)$.
\end{proof}

\begin{definition}
Suppose that $Y \subset \P^m$, $W \subset \P^n$ are two projective
varieties. The join of $Y$ and $W$, denoted as $Y\#W$, is defined
to be $\pi(C(Y)\times C(W)-0)$ where $\pi:\C^{m+n+2}-\{ 0 \}
\rightarrow \P^{m+n+1}$ is the canonical quotient map. Thus $Y\#W$
is a projective subvariety of $\P^{m+n+1}$. Define $Y \bar{\times}
W =\{ [y_0:...:y_m:w_0:...:w_n]|[y_0:...:y_m] \in Y, [w_0:...:w_n]
\in W\}=\pi(C(Y)^*\times C(X)^*)$ where $C(Y)^*=C(Y)-0$. By abuse
of notation, we will also use $Y, W$ to denote $Y \bar{\times} 0,
0 \bar{\times} W$ respectively as the images of the canonical
embedding of $Y$ and $W$ in $Y \# W$ and then we have $Y\#W=Y\cup
W\cup Y\overline{\times} W$. It is easy to see that the join of
two irreducible projective varieties is again irreducible.
\end{definition}

\begin{lemma}
Let $X \subset \P^n$, $Y \subset \P^m$ be irreducible subvarieties
and the dimension of $Y$ is $r$. Suppose that $V$ is a
$(k+r+1)$-dimensional irreducible subvariety of $Y \# X$ and $V$
intersects $X$ properly in $X\#Y$. Then $V=Y\# W$ where $W=X\cap
V$ if and only if
$$^tF_V=t^{s}F_V$$ for some $s\in \N$.
\label{chow form degree}
\end{lemma}

\begin{proof}
Suppose that $^tF_V=t^sF_V$ for some $s\in \N$. To show that
$V=Y\#W$ where $W=V\cap X$, we will first claim that $V \subset
Y\# W$. Since $V=V\cap(Y \bar{\times} X) \cup (V \cap Y) \cup
(V\cap X)$, $V\cap Y \subset Y\# W$ and $V \cap X =W$, it suffices
to show that $V \cap (Y \bar{\times} X)$ is in $Y\# W$. For every
point $a=(y_0:...:y_m:x_0:...:x_n) \in V\cap (Y \bar{\times} X)$
where $(y_0:...:y_m) \in Y$, and $(x_0:...:x_n) \in X$, we claim
that $b=(ty_0:...:ty_m:x_0:...:x_n) \in V$ for all $t \in \C$. If
$t \neq 0$, for any $k+r+2$ hyperplanes $\xi^0,...,\xi^{k+r+1}$
which contain $b$, we have
$$\sum^m_{j=0}\xi_j^ity_j+\sum_{j=0}^{n}\xi_{m+1+j}^ix_j=0$$
where $i=0,..., k+r+1$. This means that $a\in ^t\xi^{0} \cap
...\cap ^t\xi^{k+r+1}$ which implies that
$F_V(^t\underline{\xi}^{0},...,^t\underline{\xi}^{k+r+1})=0$. But
$$F_V(^t\underline{\xi}^{0},...,^t\underline{\xi}^{k+r+1})=^t
F_V(\underline{\xi}^0,...,\underline{\xi}^{k+r+1})=
t^sF_V(\underline{\xi}^0,..., \underline{\xi}^{k+r+1}) $$
therefore $F_V(\underline{\xi}^0,....,\underline{\xi}^{k+r+1})=0$.
So $b \in V$. Since those polynomials defining $V$ vanish for
infinitely many $t$, they must vanish for all t, so
$(ty_0:...:ty_m:x_0:...:x_n) \in V$ for any t. Thus
$(0:...:0:x_0:...:x_n) \in W$, therefore $b\in Y \# W$. We then
have $V \subset Y \# W$.

To show that they are the same, it is sufficient to show that $W$
is irreducible, since $V$ and $Y \# W$ have the same dimension and
both of them are irreducible. Write $W=W_1 \cup \cdots \cup W_k$
as a union of irreducible components of $W$. Since $Y\# W=
\cup^k_{i=1} Y\# W_i$ and $V \subset Y\# W$, we have $V \subset Y
\# W_i$ for some i. But then $0 \bar{\times} W_j \subset V \subset
Y\# W_i$ which implies that $W_j \subset W_i$ for all $j$. Thus
$W$ is irreducible.

Conversely, suppose that $V=Y \# W$. Fix a nonzero $t$. We are
going to show that $^tF_V$ is also a Chow form of $V$. For $a \in
V\cap \xi^0 \cap \cdots \cap \xi^{k+r+1}$, since $^{\frac{1}{t}}a
\in V \cap ^t\xi^0 \cap \cdots \cap ^t\xi^{k+r+1}$, therefore
$F_V(^t\underline{\xi}^0,...,^t\underline{\xi}^{k+r+1})=0$ which
means that
$^tF_V(\underline{\xi}^0,...,\underline{\xi}^{k+r+1})=0$.

Assume that for any $(r+1)$ hyperplanes, if $a$ is in the
intersection of them then imply $^tF_V$ evaluating at the point
defined by these $(r+1)$ hyperplanes is 0. We need to show that
$a$ is in $V$.

If $^{\frac{1}{t}}a \in \xi^0 \cap \cdots \cap \xi^{k+r+1}$, then
$a \in ^{\frac{1}{t}}\xi^0 \cap \cdots \cap
^{\frac{1}{t}}\xi^{k+r+1}$. By the assumption on $^tF_V$, we have
$^tF_V(^{\frac{1}{t}}\xi^0, ...,
^{\frac{1}{t}}\xi^{k+r+1})=F_V(\underline{\xi}^0,...,\underline{\xi}^{k+r+1})=0$.
Therefore $^ta \in V$. But $V=Y\# W$, so $a$ is in $V$. Thus
$^tF_V=cF_V$ for a nonzero constant $c$. Let
$^tF_V=F_kt^l+F_{k-1}t^{l-1}+\cdots +F_0$ then $^tF_V-cF_V\equiv
(t^l-c)F_k+(t^{l-1}-c)F_{k-1}+\cdots +(1-c)F_0\equiv 0$ for all
$t$. Because $V$ is not contained in $X$, $c$ is not $1$. So
$c=t^s$ for some $s=1,...,l$.
\end{proof}

\begin{lemma}
Let $X \subset \P^n$ be an irreducible subvariety and $Y \subset
\P^m$ be an irreducible $r$-dimensional subvariety. Suppose that
$V$ is a $(k+r+1)$-dimensional irreducible subvariety of $Y \# X$
and $V$ intersects $X$ properly in $Y\# X$. If $^tF_V=t^sF_V$ for
some $s$, then $s\geq (r+1)d$ where $d=degV$. \label{greater}
\end{lemma}

\begin{proof}
By the dimension reason, we are able to take $(r+1)$ hyperplanes
$\xi^j=\sum_{i=0}^m \xi^j_iy_i$ for $j=0,..., r$ in $\P^m$ which
do not intersect at $Y$ and $(k+1)$ hyperplanes
$\xi^j=\sum_{i=0}^n \xi^j_ix_i$ for $j=r+1,..., k+r+1$ in $\P^n$
which do not meet at $V\cap X$. So when we consider these
$(k+r+2)$ hyperplanes as hyperplanes in $\P^{n+m+1}$, they have no
common points in $V$. Therefore
$F_V(\underline{\xi}^0,...,\underline{\xi}^{k+r+1}) \neq 0$. Since
$F_V$ is a homogeneous polynomial of degree $d$ in each group of
variables where $d$ is the degree of $V$,
$F_V(\underline{\xi}^0,...,\underline{\xi}^{k+r+1})$ has a term
$g(\underline{u}^{r+1},
...,\underline{u}^{k+r+1})\prod^{r}_{j=0}\prod^m_{i=0}(\xi_i^j)^{b_{ij}}$
where $g$ is a polynomial in
$\underline{u}^{r+1}$,...,$\underline{u}^{k+r+1}$ and
$\underline{u}^j=(\xi^j_{m+1},..., \xi^j_{m+n+1})$ for $j=r+1,...,
k+r+1$. Therefore we have $\sum_{i, j} b_{ij}=(r+1)d$ which
implies that $s\geq (r+1)d$.
\end{proof}

See \cite{AN} for the following lemma:

\begin{lemma}
For a $(k+r+1)$-cycle $c$ with degree $d$ in $\P^m$,
$F_c(Au)=(detA)^dF_c(u)$ where $A=[a_{ij}]\in GL(k+r+2, \C)$ and
$Au=(\underline{u}'_0,..., \underline{u}'_{k+r+1})$ where
$u=(\underline{u}_0,..., \underline{u}_{k+r+1})$ and
$\underline{u}'_i=\sum_{j=0}^{k+r+1}a_{ij}\underline{u}_j$.
\label{det}
\end{lemma}

\begin{corollary}
Suppose that $V$ is a $(r+k+1)$-dimensional irreducible subvariety
of $\P^r\#X$ where $X\subset \P^n$. Then as a polynomial in $t$,
the degree of $^tF_V$ is less than or equal to $(r+1)d$.
\label{degree}
\end{corollary}

\begin{proof}
Denote $^t\underline{\xi}^j=(t\xi^j_0,..., t\xi^j_r,
\xi^j_{r+1},..., \xi^j_{r+n+1})$ for $j=0,..., r+k+1$. We have
$^tF_V(\underline{\xi}^0,...,
\underline{\xi}^{r+k+1})=F_V(^t\underline{\xi}^0,...,
^t\underline{\xi}^{k+r+1})$. We write $u=(^t\underline{\xi}^0,...,
^t\underline{\xi}^{r+k+1})$ in the form of a $(r+k+2)\times
(r+n+2)$-matrix
$$u=\left(%
\begin{array}{cccccc}
  t\xi^0_0 & \cdots & t\xi^0_r & \xi^0_{r+1} & \cdots & \xi^0_{r+n+1}\\
  \vdots   &        &  \vdots  & \vdots      &        & \vdots \\
  t\xi^{r+k+1}_0 & \cdots & t\xi^{r+k+1}_r & \xi^{r+k+1}_{r+1} & \cdots & \xi^{r+k+1}_{r+n+1} \\
\end{array}%
\right)$$

We can find a matrix $A\in GL(r+k+2, \C)$ such that
$$Au=\left(%
\begin{array}{cccc}
  t\xi'_0 & & 0 &  \\
  &   \ddots & &\star  \\
  0 &  & t\xi'_{r+k+1} & \\
\end{array}%
\right)$$ is in the reduced row-echelon form where $\xi'_i\in \C$
for $i=0,..., r$. Since $F_V$ is a degree $d$ homogeneous
polynomial in each group of variables, $F_V(Au)$ as a polynomial
in $t$ has degree less than or equal to $(r+1)d$. The conclusion
now follows from Lemma \ref{det}.
\end{proof}

From this Corollary and Lemma \ref{chow form degree}, we have the
following result:

\begin{theorem}
Suppose that $V$ is an irreducible subvariety of $\P^r\#X$ of
degree $d$ and $X$ is an irreducible subvariety of $\P^n$. If $V$
meets $X$ properly in $\P^r\#X$ then $V=\P^r\#W$ for $W=V\cap X$
if and only if $^tF_V=t^{(r+1)d}F_V$. \label{chow form properly}
\end{theorem}

\begin{definition}
Let $X$ be a projective variety. The suspension $\susp X$ of $X$
is defined to be $\P^0\# X$. The suspension induces a map
$\susp_*: Z_k(X) \longrightarrow Z_{k+1}(\susp X)$ from the space
of $k$-cycles on $X$ to the space of $(k+1)$-cycles on $\susp X$.
Let $\mathcal{T}^+_{k+1, d}(\susp X)$ be the open subset of
$\mathscr{C}_{k+1, d}(\susp X)$ where each $V \in
\mathcal{T}^+_{k+1, d}(\susp X)$ intersects $X \subset \susp X$
properly in $\susp X$. Let $T_{k+1}(\susp X)$ be the naive group
completion of the monoid $\mathcal{T}^+_{k+1}(\susp
X)=\coprod_{d\geq 0}\mathcal{T}^+_{k+1, d}(\susp X)$. Then
$T_{k+1}(\susp X)$ is a topological abelian subgroup of
$Z_{k+1}(\susp X)$.
\end{definition}

\begin{lemma} Suppose that $c=\sum_{i=1}^m n_iV_i\in
\mathcal{T}^+_{k+1, d}(\susp X)$. Then $c$ is in the image of
$\susp_*$ if and only if $^tF_c=t^{d}F_c$. \label{cycle chow form}
\end{lemma}

\begin{proof}
If $c$ is in the image of $\susp_*$, then each $V_i$ is in the
image of $\susp_*$, thus by Theorem \ref{chow form properly},
$^tF_{V_i}=t^{d_i}F_{V_i}$ where $d_i$ is the degree of $V_i$.
Then $$^tF_c=\Pi_{i=1}^{m}(^{t}F^{n_i}_{V_i})=\Pi_{i=1}^mt^{d_i}
F^{n_i}_{V_i}=t^{d}F_c.$$

On the other hand, write $^tF_{V_i}=t^{a_i}G_{V_i, t}$ for some
multihomogeneous polynomial $G_{V_i, t}$ which is not a multiple
of $t$. If $^tF_c=t^{d}F_c$ then
$^tF_c=\Pi^m_{i=1}t^{a_in_i}G^{n_i}_{V_i, t}= t^dF_c$. Hence
$\sum^m_{i=1} a_in_i=d=\sum^m_{i=1} n_id_i$. Since $V_i$
intersects $X$ properly in $\susp X$, by Lemma \ref{greater},
$a_i\geq d_i$ for all $i$, so we have $a_i=d_i$ for each $i$. And
this implies that $^tF_{V_i}=t^{d_i}F_{V_i}$. Therefore by Theorem
\ref{chow form properly} $c$ is in the image of $\susp_*$.
\end{proof}

Let us note that if $V\in \mathcal{T}^+_{k+1, d}(\susp X)$, we may
write $^tF_V=g_{d, V}t^{d}+\cdots +g_{1, V}t+g_{0, V}$ and $g_{d,
V}$ is not 0. Recall that $\varphi_t: \P^{n+1} \longrightarrow
\P^{n+1}$ is the map defined by $\varphi_t(x_0:...:
x_{n+1})=(tx_0:x_1:...: x_{n+1})$. In the following, we identify a
cycle $c$ with its Chow form $F_c$.

\begin{theorem}
Suppose that $X\subset \P^n$ is a projective variety. Then
$\susp_* \mathscr{C}_{k, d}(X)$ is a strong deformation retract of
$\mathcal{T}^+_{k+1, d}(\susp X)$. \label{deformation retract}
\end{theorem}

\begin{proof}
Since $^{\frac{1}{t}}F_c=t^{-d}(g_{d, c}+ g_{d-1, c}t+ \cdots +
g_{0, c}t^{d})$, $g_{d, c}+ g_{d-1, c}t+ \cdots + g_{0, c}t^{d}$
is a Chow form of $\varphi_{\frac{1}{t}}(c) \in \mathcal{T}_{k+1,
d}(\susp X)$. Define $H_d:\mathcal{T}_{k+1,d}(\susp X)\times \C
\rightarrow \mathcal{T}_{k+1,d}(\susp X)$ by
$$H_d(c,t)=g_{d,c}+\cdots +g_{1, c}t^{d-1}+g_{0,c}t^{d}$$ then $H_d$ is
a regular mapping. Since $H_d(c,1)=id$, $H_d(c,0)=g_{d, c} \in
\susp_*\mathscr{C}_{k, d}(X)$ by Lemma above where $g_{d, c}$ is
the Chow form of the cycle $\susp(c\bullet \P^n)$, $H_d(c,t)=c$
for $c \in \susp_* \mathscr{C}_{k,d}(X), \forall t \in \C$, thus
$\susp_* \mathscr{C}_{k,d}(X)$ is a strong deformation retract of
$\mathcal{T}_{k+1,d}(\susp X)$.
\end{proof}

\begin{corollary}
$\susp_* Z_k(X)$ is a deformation retract of $T_{k+1}(\susp X)$.
\label{retraction}
\end{corollary}

\begin{proof}
Let $H:\mathcal{T}^+_{k+1}(\susp X)\times
\mathcal{T}^+_{k+1}(\susp X) \longrightarrow
\mathcal{T}^+_{k+1}(\susp X)\times \mathcal{T}^+_{k+1}(\susp X)$
be defined by $H(x, y, t)=H_{d_1}(x, t)\times H_{d_2}(y, t)$ where
$d_1, d_2$ are the degree of $x, y$ respectively, so $H$ is
continuous. We observe that since $F_{x+y}=F_x\cdot F_y$, we have
$H_d(x+y, t)=H_{d_1}(x, t)H_{d_2}(y, t)$ where $d=d_1+d_2$. If
$(x, y) \sim (w, z)$, i.e. $x+z=w+y$, then we have the equality
$H_k(x+z, t)=H_k(w+y, t)$ where $k=d_1+d_2$. Thus $H_{d_1}(x,
t)\times H_{e_2}(z, t)=H_{e_1}(w, t)\times H_{d_2}(y, t)$ where
$e_1, e_2$ are the degree of $w,z$ respectively. This means that
$(H_{d_1}(x, t), H_{d_2}(y, t))\sim (H_{e_1}(w, t), H_{e_2}(z,
t))$. Therefore, $H$ reduces to a map
$$H:T_{k+1}(\susp X)\times \C \longrightarrow T_{k+1}(\susp X)$$
with the properties: $H(c, 1)=c$; $H(c, 0)\in \susp_* Z_k(X)$;
$H(c, t)=c$ for all $c \in \susp_* Z_k(X)$. Thus $\susp_* Z_k(X)$
is a strong deformation retract of $T_{k+1}(\susp X)$.
\end{proof}

\section{The Holomorphic Taffy for totally real cycles}
Suppose that $X\subset \RP^n$ is a totally real projective
variety. Let $RT_{p+1}(\susp X)=T_{p+1}(\susp X_{\C})\cap
RZ_{p+1}(\susp X)$, $\overline{RT}_{p+1}(\susp X)$ be the closure
of $RT_{p+1}(\susp X)$ in $T_{p+1}(\susp X_{\C})$, and
$\overline{RT}_{p+1, d}(\susp X)=T_{p+1, d}(\susp X_{\C})\cap
\overline{RT}_{p+1, d}(\susp X)$. Let $\overline{RT}^+_{p+1,
d}(\susp X)$ be the collection of all effective cycles in
$\overline{RT}_{p+1, d}(\susp X)$ and define
$$S\mathscr{C}_{p, d}(X)=\{c\bullet \P^n| c\in \overline{RT}^+_{p+1,
d}(\susp X)\}$$ where $\bullet$ is the intersection product.

Let $SZ_p(X)$ be the naive group completion of the monoid
$S\mathscr{C}_p(X)=\coprod_{d\geq 0}S\mathscr{C}_{p, d}(X).$ Then
$RZ_p(X)$ is a subgroup of $SZ_p(X)$.

Let $H_d:\mathcal{T}_{p+1,d}(\susp X)\times \C \rightarrow
\mathcal{T}_{p+1,d}(\susp X)$ be the deformation retraction
defined in Theorem \ref{deformation retract}. We have the
following result:

\begin{theorem}
The restriction of $H_d$ to $\overline{RT}^+_{p+1, d}(\susp
X)\times \R$ is a strong deformation retraction of
$\overline{RT}^+_{p+1, d}(\susp X)$ onto $\susp_* S\mathscr{C}_{p,
d}(X)$. \label{real taffy}
\end{theorem}

\begin{proof}
For $c\in RT^+_{p+1, d}(\susp X)$, $\varphi_t(c)\in RT_{p+1,
d}(\susp X)$ for $t\in \R\backslash \{0\}$. So $\varphi_0(c)\in
\overline{RT}_{p+1, d}(\susp X)$. Therefore $\varphi_t(c)\in
\overline{RT}_{p+1, d}(\susp X)$ for $c\in \overline{RT}_{p+1,
d}(\susp X)$. We have $H_d(c, 0)=g_{d, c}$ which is the Chow form
of $\susp(c\bullet \P^n)$ for $c\in \overline{RT}_{p+1, d}(\susp
X)$, and by definition $H_d(c, 0)$ is in $\susp_* S\mathscr{C}_{p,
d}(X)$.
\end{proof}

\begin{corollary}
$\susp_* SZ_p(X)$ is a strong deformation retract of
$\overline{RT}_{p+1}(\susp X)$. \label{retract}
\end{corollary}

\begin{proof}
From the similar analysis as in the proof of Theorem \ref{real
taffy}, the restriction of the map $H$ in Corollary
\ref{retraction} is a strong deformation retraction of
$\overline{RT}_{p+1}(\susp X)$ onto $\susp_* SZ_p(X)$.
\end{proof}

\section{Critical polynomials}
\begin{definition}
Suppose that $f \in \R[x]$ is a real polynomial of degree n. We
define a sequence of real polynomials $(f_0, f_1,..., f_k)$ as
following: let $f_0=f$, $f_1=f'$ be the derivative of $f$, then by
using the Euclidean algorithm, we can find $q_i, f_i \in \R[x]$
such that $f_0=q_1f_1-f_2$, $f_1=q_2f_2-f_3$, ....,
$f_{k-1}=q_kf_k-0$ and $f_{k-1} \neq 0$ where $deg f_{i+1} < deg
f_i$ for all $i=1,...,k-1$. The sequence $(f_0,f_1,...,f_k)$ is
called the Sturm sequence of $f$.
\end{definition}

Given a sequence of real numbers $s=(a_1,...,a_n)$, if each $a_i$
is not zero, we define the number of sign changes of $s$ to be the
number $sgn(s)=|\{a_i|a_ia_{i+1}<0\}|$. If some $a_i$ are zeros,
deleting all zeros, we get a new sequence $s'$ and define
$sgn(s)=sgn(s')$.

The following Sturm theorem can be found in \cite{BCR}, Corollary
1.2.10.
\begin{theorem}
Suppose that $f\in \R[x]$ is a real polynomial of degree n and let
$s=(f_0,...,f_k)$ be the Sturm sequence of $f$. For $a \in \R$,
define $sgnf_a=sgn(f_0(a),..., f_k(a))$. Then the number of
distinct real roots of $f$ in the interval $(a, b) \subset \R$ is
equal to $sgnf_a-sgnf_b$.
\end{theorem}

\begin{corollary}
The number of distinct real roots of a real polynomial is $\lim_{a
\to \infty}(sgnf_{-a}-sgnf_a)$.
\end{corollary}

\begin{definition}
For a monomial $a^{i_1}_1...a^{i_n}_n \in \R[a_1,...,a_n]$, we
define its substituted homogeneous degree (shd) to be
$$shd(a^{i_1}_1...a^{i_n}_n)=\sum^n_{j=1}j\cdot i_j \in \N.$$ A
polynomial $p_i(a_1,..., a_n)=\sum
c_{i_1,...,i_n}a^{i_1}_1...a^{i_n}_n \in \R[a_1,..., a_n]$ is
called a substitutable homogeneous polynomial if all the terms
$a^{i_1}_1...a^{i_n}_n$ have the same substituted homogeneous
degree. For a rational function $\frac{p(a_1,..., a_n)}{q(a_1,...,
a_n)}$, we say that it is a substitutable homogeneous rational
function if there exist two substitutable homogeneous polynomials
$p_1$, $q_1$ in $\R[a_1,..., a_n]$ such that
$\frac{p}{q}=\frac{p_1}{q_1}$. We define the substituted
homogeneous polynomial degree of $\frac{p}{q}$ to be
$shd(\frac{p}{q}) =shd(p_1)-shd(q_1)$.
\end{definition}

If $\frac{p_1}{q_1}=\frac{p_2}{q_2}$ where $p_1$, $q_1$, $p_2$,
$q_2$ are all substitutable homogeneous polynomials then
$shd(p_1)+shd(q_2)=shd(p_2)+shd(q_1)$. Therefore
$shd(p_1)-shd(q_1)=shd(p_2)-shd(q_2)$. So we know that the $shd$
for a substitutable homogenous rational functions is well-defined.

The following is a straight forward calculation.
\begin{proposition}
Suppose that $f(a_1,..., a_n)$ is a substitutable homogenous
polynomial, and $g_i(x_1,.., x_m)$ is a homogeneous polynomial of
degree $i$ for $i=1,..., n$, then
$$f(g_1(x_1,..., x_m),..., g_n(x_1,..., x_m))$$ is a homogeneous polynomial in
$x_1,..., x_m$ with degree $shd(f)$.
\end{proposition}

\begin{definition}
Let $A(x)= p_0(a_1,..., a_n)x^d+ \cdots + p_d(a_1,..., a_n)$,
$B(x)=q_0(a_1,..., a_n)x^{d-1}+ \cdots + q_{d-1}(a_1,..., a_n)$
where $p_i(a_1,..., a_n)$, $q_i(a_1,..., a_n)$ are all
substitutable homogeneous rational functions for $i=0,..., d$ and
$j=0,..., d-1$. If

1. $shdq_i-shdp_i=c$ a constant, for $i=0,..., d-1$

2. $shdp_i=i+shdp_0$, for $i=0,..., d$

then $(A, B)$ is called a substitutable pair.
\end{definition}

\begin{proposition}
Let $A(x)$, $B(x)$ as above and $(A, B)$ a substitutable pair. By
doing the Euclidean algorithm, we may express $A(x)=B(x)D(x)+C(x)$
where $degC(x)<degB(x)$. Then $(B, C)$ is also a substitutable
pair.
\end{proposition}

\begin{proof}
Let
$r_{i-2}=\frac{p_iq_0^2-p_0q_0q_i-p_1q_0q_{i-1}+p_0q_1q_{i-1}}{q_0^2}$
for $i=2,..., d-1$ and
$r_{d-2}=\frac{p_dq_0^2-p_1q_0q_{d-1}+p_0q_1q_{d-1}}{q_0^2}$. From
the Euclidean algorithm, we get $C(x)=r_0(a_1,..., a_n)x^{d-2}+
\cdots+r_d(a_1,..., a_n)$. And it is easy to check that
$r_i(a_1,..., a_n)$ is a substitutable homogeneous rational
function with $shd(r_i)=shd(p_{i+2})$ for $i=0,..., d-2$. Since
$shd(r_i)-shd(q_i) =shd(p_{i+2})-shd(q_i)=2-C$, for all $i=0,..,
d-2$ and $shd(q_i)=i+shd(q_0)$, for all $i=0,..., d-1$, thus $(B,
C)$ is a substitutable pair.
\end{proof}

With all the work above, we can show that any consecutive term of
the Sturm sequence of a polynomial is a substitutable pair.

\begin{proposition}
Suppose that $f(x)=x^d+a_1x^{d-1}+ \cdots +a_d$ is a real
polynomial and $(f_0,f_1,..., f_d)$ is the Sturm sequence of $f$.
Consider $f_i$ as a polynomial in $a_1,..., a_d$ for all $i$. Then
$(f_i, f_{i+1})$ is a substitutable pair. \label{codim}
\end{proposition}

\begin{theorem}
Suppose that $f(x)=x^d+a_1x^{d-1}+ \cdots +a_d$ is a real
polynomial. Then there exist substitutable polynomials $F_2,...,
F_d$ such that $f$ has d distinct real roots if and only if
$F_2(a_1,..., a_d)>0, \cdots, F_d(a_1,..., a_d)>0$.
\label{critical poly}
\end{theorem}

\begin{proof}
Suppose that $(f_0,f_1,..., f_d)$ is the Sturm sequence of $f$. By
Proposition \ref{codim}, the leading coefficient of $f_i$ is a
substitutable rational function which can be expressed as
$\frac{F_i(a_1,..., a_n)}{(w_i(a_1,..., a_n))^2}$ where $F_i$ and
$w_i$ are substitutable homogeneous polynomials. Considering the
cases of $d$ being even and odd separately by taking $x$ to be very
negative, the Sturm theorem says that $f$ has $d$ distinct real
roots if and only if
$$\left\{
    \begin{array}{rl}

F_2(a_1,..., a_n)>0   \\

\vdots \\
F_d(a_1,..., a_n)>0
\end{array}\right.$$

\end{proof}
We call $F_2,..., F_d$ the critical polynomials of degree $d$
polynomials. This theorem will be used in proving the openness of
some sets in \emph{magic fans}.

For a real polynomial $f(x)\in \R[x]$ of degree $d$ with leading
coefficient 1 , we identify it with a point in $\R^d$ by taking its
coefficients as coordinates. Then the above result gives us the
following byproduct.

\begin{corollary}
The set $S_n=\{ (a_1,..., a_n) \in \R^n| x^n+a_1x^{n-1}+...+a_n=0
\mbox{ has n distinct real roots }\}$ is an open semi-algebraic
set of $\R^n$.
\end{corollary}

\section{Magic Fan}
\begin{definition}
A homogeneous polynomial $f(x_1,..., x_n)$ is said to be positive
if $f(x_1,..., x_n)>0$ for all $(x_1,..., x_n) \in \R^n-\{0\}$.
\end{definition}

Let us designate the set of all effective divisors of degree $d$
in $\P^n$ defined by real homogenous polynomials by $Div_{n, d}$.
By using Veronese embedding, we can identify $Div_{n, d}$ with
$\RP^{\binom{n+d}{d}-1}$ and topologize it with the standard
topology of $\RP^{\binom{n+d}{d}-1}$. We often abuse a homogenous
polynomial with the effective divisor defined by it in a
projective space.

\begin{definition}
Let $Div'_{n,d}$ be the collection of effective divisors of degree
$d$ of $\P^{n+1}$ defined by real homogeneous polynomials in the
form
$$G(x_0,...,
x_n)=x_0^d+\sum_{\overset{i_0+...+i_n=d}{i_0<d}}c_{i_0,...,i_n}x^{i_0}_0\cdots
x^{i_n}_n$$ where $c_{i_0,..., i_n} \in \R.$ Since $Div'_{n, d}$
is the set of all divisors which do not contain $[1:0:\cdots:0]$,
$Div'_{n, d}$ is an open subset of $Div_{n, d}$.
\end{definition}

\begin{definition}
Let $E_{n,d}$ be the subset of $Div'_{n,d}$ consisting of all
$G(x_0,..., x_n)$ such that for any $(x_1,..., x_n) \in
\R^n-\{0\}$, as a polynomial in $x_0$, $G(x_0,..., x_n)$ has d
distinct real roots.
\end{definition}

\begin{example}
It is obvious that $G(x_0,..., x_n)=(x^2_0-(x^2_1+...+
x^2_n))\cdots (x^2_0-d(x^2_1+...+x^2_n))$ is in $E_{n,2d}$ and
$H(x_0,..., x_n)=x_0G(x_0,..., x_n)$ is in $E_{n, 2d+1}$. Thus
$E_{n, d}$ is not empty.
\end{example}

Our goal here is to show that $E_{n,d}$ is an open subset of
$Div_{n, d}$. It suffices to show that $E_{n, d}$ is an open
subset of $Div'_{n, d}$. For any $G$ in $E_{n,d}$, we need to find
an $\epsilon>0$ such that if $|c'_{i_0,..., i_n}-c_{i_0,...,
i_n}|< \epsilon$ for any $i_0, ..., i_n$ where $c'_{i_0,...,
i_n}$, $c_{i_0,...,i_n}$ are coefficients of the term
$x^{i_0}_0\cdots x^{i_n}_n$ of $G'$ and $G$ respectively, then
$G'$ is in $E_{n, d}$ .

\begin{lemma}
Suppose that $H(x_1,..., x_n)=\sum_{i_1+\cdots +
i_n=k}h_{i_1,...,i_n}x^{i_1}_1\cdots x^{i_n}_n \in \R[x_1,...,
x_n]$ is a positive homogeneous polynomial, then there is an
$\epsilon>0$ such that if $|h'_{i_1,..., i_n}-h_{i_1,...,
i_n}|<\epsilon$, for all $i_1,.., i_n$, the homogenous polynomial
$H'(x_1,..., x_n)=\sum_{i_1+\cdots +i_n=k}h'_{i_1,...,
i_n}x^{i_1}_1\cdots x^{i_n}_n$ is also positive.
\end{lemma}

\begin{proof}
Consider $H|_{S^{n-1}}$, the restriction of $H$ to the unit sphere
in $\R^n$. Let $\delta$ be the minimum value of $H|_{S^{n-1}}$,
thus $\delta>0$. Let $M=\binom{n+k-1}{k}$ and $\epsilon =
\frac{\frac{1}{2}\delta}{M}$. Then for any $|h'_{i_1,...,
i_n}-h_{i_1,..., i_n}|<\epsilon$ where $H'(x_1,...,
x_n)=\sum_{i_1+\cdots +i_n=k}h'_{i_1,..., i_n}x^{i_1}_1\cdots
x^{i_n}_n$, we have
$|H'|_{S^{n-1}}(X)-H|_{S^{n-1}}(X)|<\frac{1}{2}\delta$ for $X \in
S^{n-1}$. Thus $H'|_{S^{n-1}}(X)>\frac{1}{2}\delta$ for $X \in
S^{n-1}$. Then for any $X\in\R^n-\{0\}$,
$H'(X)=|X|^kH'|_{S^{n-1}}(\frac{X}{|X|})>|X|^k\frac{\delta}{2}>0$
since $H'$ is homogeneous of degree $k$. Thus $H'$ is a positive
homogeneous polynomial.
\end{proof}

\begin{lemma}
\label{opensubset} $E_{n,d}$ is an open subset of $Div'_{n,d}$.
\end{lemma}

\begin{proof}
Let $F_j$, $j=2,..., d$ be the critical polynomials of degree $d$
polynomials and $k_j$ be the substituted homogeneous degree of
$F_j$. Suppose that
$$G(x_0,..., x_n)=x^d_0+p_1(x_1,..., x_n)x^{d-1}_0+\cdots
p_d(x_1,..., x_n)$$
$$=x^d_0+\sum_{\overset{i_0+\cdots +i_n =d}{i_0<d}}c_{i_0,..., i_n}
x^{i_0}_0\cdots x^{i_n}_n \in E_{n,d}.$$ Then $F_j(p_1(x_1,...,
x_n),..., p_d(x_1,..., x_n)) =\sum_{r_1+\cdots
+r_n=k_j}t_{r_1,..., r_n}x^{r_1}_1\cdots x^{r_n}_n$ is a positive
homogeneous polynomial. By Lemma above, there is an $\epsilon>0$
such that if $|t'_{r_1,..., r_n}-t_{r_1,..., r_n}|<\epsilon$, then
$F'(x_1,..., x_n)=\sum_{r_1+ \cdots +r_n=k_j}t'_{r_1,...,
r_n}x^{r_1}_1\cdots x^{r_n}_n$ is a positive homogenous
polynomial. Let
$$G'(x_0,..., x_n) =x^d_0+q_1(x_1,..., x_n)x^{d-1}_0+\cdots
+q_d(x_1,..., x_n)$$
$$=x^d_0+\sum_{\overset{i_0+\cdots +i_n =d}{i_0<d}}c'_
{i_0,..., i_n}x^{i_0}_0\cdots x^{i_n}_n \in Div'_{n,d}$$
Considering $G'$ as a polynomial in $x_0$ and write the critical
polynomial
$$F_j(q_1(x_1,..., x_n),..., q_d(x_1,..., x_n))= \sum_{r_1+\cdots
+r_n=k_j}t'_{r_1,..., r_n}x^{r_1}_1\cdots x^{r_n}_n.$$ Then
$t'_{r_1,..., r_n}$ is a polynomial in variables $\{c'_{i_0,...,
i_n}|i_0+ \cdots +i_n=k_j, i_0<d\}$. Thus there is a $\delta>0$
such that if $|c'_{i_0,..., i_n}-c_{i_0,..., i_n}|<\delta$, then
$|t'_{i_0,..., i_n}-t_{i_0,..., i_n}|<\epsilon$ and by above
assumption
$$F_i(q_1(x_1,..., x_n),..., q_d(x_1,..., x_n))>0$$ for all $(x_1,...,
x_n) \in \R^n-\{0\}$. Therefore by Theorem \ref{critical poly},
$G'(x_0,..., x_n)$ has d distinct real roots for all $(x_1,...,
x_n) \in \R^n-\{0\}$. Thus $E_{n,d}$ is open in $Div'_{n,d}$.
\end{proof}

For
$$f(x_0,..., x_n)=x^d_0+\sum_{\overset{i_0+\cdots+i_n=d}{i_0<d}}
a_{i_0,...,i_n}x^{i_0}_0\cdots x^{i_n}_n \in Div'_{n, d},$$ let
$$f_t(x_0,..., x_n)=t^df(\frac{1}{t}x_0,x_1,..., x_n)$$ and denote
$tD$ to be the divisor in $\P^n$ determined by $f_t$ where $D$ is
the divisor determined by $f$. We use $0D$ to denote the cycle
$deg(D)\cdot \P^{n-1}$, the hyperplane which has the same degree
as $D$'s, and from the limiting process as $t$ approaches 0, we
denote $f_0=x_0^d$.

\begin{definition}
Let $Div''_{n, d}$ be the subset of $E_{n, d}$ consisting all
homogeneous polynomial $f$ of degree $d$ such that $f_t(1, 1,
0,..., 0)\neq 0$ for all $t\in (0,1]$.
\end{definition}

\begin{proposition}
$Div''_{n, d}$ is an open subset of $Div_{n, d}$. \label{open
subset}
\end{proposition}

\begin{proof}
Let $f(x_0,...,
x_n)=x^d_0+\sum_{\overset{i_0+\cdots+i_n=d}{i_0<d}}
c_{i_0,...,i_n}x^{i_0}_0\cdots x^{i_n}_n$ and $g(x_0,..., x_n)=\\
x^d_0+\sum_{\overset{i_0+\cdots+i_n=d}{i_0<d}}
c'_{i_0,...,i_n}x^{i_0}_0\cdots x^{i_n}_n$ be two polynomials in
$E_{n, d}$. Then $f_t(1, 1, 0,...,
0)=1+\sum_{\overset{i_0+i_1=d}{i_0<d}} c_{i_0,...,i_n}t^{i_1}$ and
$g_t(1, 1, 0,..., 0)=1+\sum_{\overset{i_0+i_1=d}{i_0<d}}
c_{i_0,...,i_n}t^{i_1}$. If $f\in Div''_{n, d}$, then
$$|f_t(1,1,0,..., 0)-g_t(1,1,0,..., 0)|\leq
\sum_{\overset{i_0+i_1=d}{i_0<d}}|c_{i_0,...,i_n}-c'_{i_0,...,i_n}|$$
for $t\in (0, 1]$. Therefore, if $g$ is close enough to $f$, $g$
is in $Div''_{n, d}$.
\end{proof}

\begin{definition}
Let $c$ be a real $r$-cycle on $\P^n$ and let
$\mathcal{O}_c(d)_{\R}$ be the set of polynomials in $Div_{n, d}$
restricted to the support of $c$. For example,
$\mathcal{O}_{\P^r}(d)_{\R}$ is the set of all real homogeneous
polynomials over $\P^r$ with degree less than or equal to $d$.
Therefore,
$dim_{\R}\mathcal{O}_{\P^r}(d)_{\R}=\binom{r+d+1}{d}-1$.

Let $x_{\infty}=(1:0:...:0), x_{11}=(1:1:0:...:0)$ be two points
in $\P^{n+1}$. For an effective real $r$-cycle $c$ on $X\subset
\P^n$, let $\beta_{n+1, d}(c)$ be the subset of $Div_{n+1, d}$
consisting of real divisors $D$ such that $x_{11}\#\P^n$ contains
an irreducible component of $(x_{\infty}\#c)\cap D$, or
equivalently, $D$ contains an irreducible component of
$(x_{11}\#\P^n)\cap (x_{\infty}\#c)$.
\end{definition}

\begin{lemma}
For a real cycle $c$ on $\P^n$,  the real codimension
$codim_{\R}(\beta_{n+1, d}(c))$ of $\beta_{n+1, d}(c)$ is greater
or equal to $\binom{r+d+1}{d}-1$. \label{codimension}
\end{lemma}

\begin{proof}
Let $V$ be an irreducible component of $(x_{11}\#\P^n)\cap
(x_{\infty}\#c)$. Take a real linear subspace $L$ of dimension
$n+1-r$ in $\P^{n+1}$ which does not intersect with $V$. We may
assume that the projection defined by $L$ is $\pi: \P^{n+1}-L
\longrightarrow \P^r$ and thus $\pi(V)=\P^r$. If $f$ is a real
homogenous polynomial of degree $d$ defined over $\P^{n+1}$ and
$f|_{\P^r}$ is not the zero polynomial, then we may define
$g(x_0,..., x_{n+1})=f(x_0,..., x_r, 0,..., 0)$ which is not a
zero polynomial on $V$. Hence the divisor defined by $g$ does not
contain $V$. Therefore, $dim_{\R}\mathcal{O}_V(d)_{\R}\geq
dim_{\R}\mathcal{O}_{\P^r}(d)_{\R}\geq \binom{r+d+1}{d}-1$.

Let $(x_{11}\#\P^n)\cap (x_{\infty}\#c)=\sum^m_{i=1}n_iV_i$, and
consider the restriction map $\psi_i:Div_{n+1, d} \longrightarrow
\mathcal{O}_{V_i}(d)_{\R}$. Then $codim_{\R}\beta_{n+1,
d}(c)=\underset{i}{min}\{ dim Im \Psi_i\}=\underset{i}{min}\{ dim
\mathcal{O}_{V_i}(d)_{\R}\}\geq \binom{r+d+1}{d}-1$.
\end{proof}

\begin{definition}
Suppose that $c$ is a real $r$-cycle on $\P^{n+1}$. Define
$\beta''_{n+1, d}(c)=\beta_{n+1, d}(c)\cap Div''_{n+1, d}$.
\end{definition}

\begin{corollary}
$codim_{\R}\beta''_{n+1, d}(c) \geq \binom{r+d+1}{d}-1$.
\end{corollary}

\begin{proof}
By Lemma \ref{open subset}, $Div''_{n+1, d}$ is an open subset of
$Div_{n+1, d}$, thus we have $codim_{\R}\beta''_{n+1,d}(c)\geq
codim_{\R}\beta_{n+1, d}(c)\geq \binom{r+d+1}{d}-1$.
\end{proof}

For a real projective variety $V\subset \P^n$, we use
$Re(V)\subset \RP^n$ to denote the set of real points of $V$.

\begin{lemma}
Let $f_j(x_0,..., x_n)=\sum_{i=0}^n a_ix_i \in \R[x_0,..., x_n]$,
for $j=0,..., m$ and define a projection $\pi:\P^n-L
\longrightarrow \P^m$ by $\pi=(f_0,..., f_m)$ where $L$ is the
linear subspace defined by $f_0,..., f_m$. If $V \subset \P^n-L$
is an irreducible projective variety of dimension $p$ and
$dim_{\R}Re(V)=p$, then $dim_{\R}Re(\pi(V))=p$. In other words, if
$V$ is the complexification of $Re(V)$, then $\pi(V)$ is the
complexification of $Re(\pi(V))$.
\end{lemma}

\begin{proof}
Assume that $dim_{\R}Re(\pi(V))<p$. Let $W$ be the
complexification of $Re(\pi(V))$ then $W \subsetneq \pi(V)$. Each
component of $\pi$ is a real linear function, thus $\pi(Re(V))
\subset Re(\pi(V))\subset W$. From here we have $Re(V)\subset
\pi^{-1}(W)\cap V$. Since $W$ is a proper subset of $\pi(V)$, thus
$\pi^{-1}(W)\cap V \subsetneq V$. But $\pi^{-1}(W)\cap V$ is a
proper subvariety of $V$ and contains $Re(V)$, so $V$ is not the
complexification of $Re(V)$. This contradicts to the hypothesis.
\end{proof}

Suppose that $X\subset \P^n$ is a real projective variety. For
$t\in \C$, $c\in Z_{r+1}(\susp X)$ and $D\in Div''_{n+2,
d}-\beta''_{n+2, d}(c)$, define $\Psi_{tD}(c)=\pi_{1*}(\susp c
\bullet tD)$ where $\pi_1:\P^{n+2}-x_{11} \longrightarrow
\P^{n+2}$ is the projection with center $x_{11}$ and $\susp$ is
the map induced by joining with $x_{\infty}$. Let $\pi_{\infty}:
\P^{n+2}-x_{\infty} \longrightarrow \P^{n+1}$ be the projection
with center $x_{\infty}$.

\begin{lemma}
Suppose that $X\subset \P^n$ is a real projective variety and
$c\in RZ_{r+1, e}(\susp X)$. For $D\in Div''_{n+2,
d}-\beta''_{n+2, d}(c)$, $\Psi_{tD}(c) \in RZ_{r+1, de}(\susp X)$
for $ t\in [0, 1]$. \label{real map}
\end{lemma}

\begin{proof}
Let $V$ be an irreducible component of $(\susp c) \cap tD$ where
$t\in [0, 1]$ and $Z=\pi_{\infty}(V)$. So $Z$ is an irreducible
component of the support of $c$ and therefore $dim_{\R}Re(Z)=r+1$.
The projection $\pi_{\infty}$ is a real map hence $\pi_{\infty}(Re
(V)) \subset Re(Z)$. Since $D\in Div''_{n+2, d}$, for any real
point $q \in Re(Z)$, $(x_{\infty}\#q) \cap V$ are some real
points, thus $\pi_{\infty}(Re(V))=Re(Z)$. But
$$r+1 \geq dim_{\R}Re(V) \geq dim_{\R}\pi_{\infty}(Re(V))=dim_{\R}Re(Z)=r+1$$
thus $V$ is the complexification of $Re(V)$ and by the Lemma
above, $\pi_1(V)$ is the complexification of $Re(\pi_1(V))$.
\end{proof}

\begin{proposition}
Let $X\subset \RP^n$ be a totally real projective variety and $K
\subset \overline{RZ}_{r+1, e}(\susp X)$ be a compact subset.
Define $\beta_{n+2, d}(K)=\cup_{Z \in K}\beta_{n+2, d}(Z)$. If
$\binom{r+d+1}{d}>dim_{\R}(K)+1$, then there exists $D\in
Div''_{n+2, d}$ and a continuous map
$$\Psi_{tD}:K \longrightarrow
\overline{RZ}_{r+1, de}(\susp X) \mbox{ for } t\in[0, 1],$$
defined by
$$\Psi_{tD}(Z)=\pi_{1*}((\susp Z)\bullet tD)$$ with the
following properties:
\begin{enumerate}
    \item $\Psi_{0D}(Z)=d\cdot Z$
    \item $\Psi_{tD}(K) \subset \overline{RT}_{r+1,
de}(\susp X) \mbox{ for } t\in (0,1].$
\end{enumerate}
\label{divisor map}
\end{proposition}

\begin{proof}
From Lemma \ref{codimension} and the assumption, we have
$codim_{\R}\beta''_{n+2, d}(K)= max\{codim_{\R}\beta''_{n+2,
d}(Z)|Z\in K\}-dim_{\R}(K)>\binom{r+d+1}{d}-1-dim_{\R}(K)>0$.
Therefore we can find a $D\in Div''_{n+2, d}$ which is not in
$\beta''_{n+2, d}(K)$. From Proposition 3.5 in \cite{F}, we know
that the map
$$\Psi_{tD}:K \longrightarrow
Z_{r+1, de}(\susp X)_{\R},$$ defined by
$$\Psi_{tD}(Z)=\pi_1((\susp Z\bullet tD)$$ is continuous and
$\Psi_{tD}(Z)$ meets $X$ properly in $\susp X$. Now let us show
that the image of $\Psi$ is actually contained in
$\overline{RZ}_{r+1, de}(\susp X)$. By Lemma \ref{real map}, for
any $Z\in K\cap RZ_{r+1, e}(\susp X)$, $\Psi_{tD}(Z)\in RT_{r+1,
de}(\susp X)$ for $t\in (0, 1]$. Therefore, if $Z\in
\overline{RZ}_{r+1, e}(\susp X)$, $\Psi_{tD}(Z)\in
\overline{RT}_{r+1, de}(Z)$ for $t\in (0, 1]$.
\end{proof}

Recall that the topology of $Z_{r+1}(\susp X_{\C})_{\R}$ is
defined from a compactly filtered filtration (see \cite{T1}, page
8)
$$K_{r+1, 1}(\susp X_{\C})_{\R} \subset K_{r+1, 2}(\susp X_{\C})_{\R}
\subset \cdots =Z_{r+1}(\susp X_{\C})_{\R}$$ by the weak topology.
And by Lemma 2.2 of \cite{T1}, this filtration is locally compact
which means that for a compact set $K\subset Z_{r+1}(\susp
X_{\C})_{\R}$, there is a number $k$ such that $K \subset K_{r+1,
k}(\susp X_{\C})$. Therefore the filtration formed by $RK_{r+1,
d}(\susp X)=K_{r+1, d}(\susp X_{\C})_{\R}\cap
\overline{RZ}_{r+1}(\susp X)$ of $\overline{RZ}_{r+1}(\susp X)$ is
locally compact.

\begin{theorem}
The inclusion map $i: \overline{RT}_{r+1}(\susp X) \longrightarrow
\overline{RZ}_{r+1}(\susp X)$ is a weak homotopy equivalence.
\end{theorem}

\begin{proof}
For surjectivity, consider a continuous map $f:S^m \longrightarrow
\overline{RZ}_{r+1}(\susp X)$, then $f(S^m)\subset
RK_{r+1,e}(\susp X)$ for some $e$. By Proposition \ref{divisor
map}, for $d$ large enough, there is a map
$$\Psi_{tD}:K \longrightarrow
\overline{RZ}_{r+1, de}(\susp X) \mbox{ for } t\in[0, 1]$$ which
is a homotopy between the map $d\cdot f$ and $g_d=\Psi_{1D}(f)$
where $g_d:S^m \longrightarrow \overline{RT}_{r+1, de}(\susp X)$.
Thus $i_*([g_{d+1}]-[g_d])=[f]$. For injectivity, consider a
continuous map of pair $f:(D^{m+1}, S^m)\longrightarrow
(\overline{RZ}_{r+1}(\susp X), \overline{RT}_{r+1}(\susp X))$
which is nullhomotopic. Take $d$ large enough and then by
Proposition \ref{divisor map}, we have a map
$\Psi_{tD}:(f(D^{m+1}), f(S^m))\longrightarrow
(\overline{RZ}_{r+1}(\susp X), \overline{RT}_{r+1}(\susp X))$
which is a homotopy between $\Psi_{0D}(f)=d\cdot f$ and
$g_d=\Psi_{1D}(f)$ where $g_d:(D^{m+1}, S^m)\longrightarrow
(\overline{RT}_{r+1}(\susp X), \overline{RT}_{r+1}(\susp X))$. So
$g_d$ is nullhomotopic. Thus $i_*(g_d)=0$ and
$i_*([g_{d+1}]-[g_d])=[f]=0$.
\end{proof}

Combining this Theorem with Corollary \ref{retract}, we get
\begin{theorem}
Suppose that $X\subset \RP^n$ is a totally real projective
variety. The suspension map $\susp_*: SZ_r(X) \longrightarrow
\overline{RZ}_{r+1}(\susp X)$ induces a weak homotopy equivalence.
\end{theorem}

\begin{example}
$$\overline{RZ}_1(\RP^{n+1}) \mbox{ is weak homotopy equivalent to }
\prod_{i=0}^n\prod_{j=0}^iK(I_{i, j}, i+j)$$ where
$$I_{i, j}=
\left\{%
\begin{array}{ll}
    0, & \hbox{if $j$ is odd or $j>i$;} \\
    \Z, & \hbox{if $j=i$ and $j$ is even ;} \\
    \Z_2, & \hbox{if $j<i$ and $j$ is even.} \\
\end{array}%
\right.
$$
and $K(I_{i, j}, i+j)$ is the Eilenberg-Mac Lane space.
\end{example}
\begin{proof}
For any real cycle $c\in Z_0(\P^n)_{\R}$, we can take a cycle $c'
\in RT_1(\P^{n+1})$ such that $c=c'\bullet \P^n$. Therefore
$SZ_0(\RP^n)=Z_0(\P^n)_{\R}$. With the Theorem above, the
conclusion now follows from Theorem 3.3 in \cite{LLM}.
\end{proof}

\begin{acknowledgement}
The author thanks Blaine Lawson for his patience in listening to the
proof and his guidance for the author's graduate study in Stony
Brook. He also thanks the referee for his/her comments.
\end{acknowledgement}

\bibliographystyle{amsplain}

\begin{thebibliography}{1}
\bibitem{AK} S. Akbulut, H. King, The Topology of Real Algebraic
Sets, L'Enseignement Math., 29 (1983), 221-261.

\bibitem{AN} A. Andreotti, F. Norguet, La convexit\'{e}
holomorphe dans l'espace analytique des cycles d'une
vari\'{e}t\'{e} alg\'{e}brique, Ann. Scuola Norm. Sup. Pisa, 21
(1967), 31-82.

\bibitem{BCR} J. Bochnak, M. Coste and M.-F. Roy, Real Algebraic
Geometry, Ergebnisse der Math. und ihrer Grenzgeb. Folge 3, Vol.
36, Berlin Heidelberg, New York, Springer, 1998.

\bibitem{Ben} M. Ben-Or, Lower bounds for algebraic computation
tress, Proc. 15th ACM Annual Symp. on Theory of Comput., 80-86
(1983).

\bibitem{F} E. Friedlander, Algebraic cycles, Chow
varieties, and Lawson homology, Compositio Math. 77 (1991), 55-93.

\bibitem{L} H.B. Lawson, Algebraic cycles and homotopy
theory, Annals of Math. 129 (1989), 253-291.

\bibitem{LLM} H.B. Lawson, P. Lima-Filho, M. Michelsohn, Algebraic
cycles and the classical groups I: real cyclces, Topology, {\bf
42}, (2003) 467-506.

\bibitem{M} J. Milnor,  On the Betti numbers of real
varieties, Proc. Amer. Math. Soc. 15 (1964), 275-280.

\bibitem{S} P. Samuel, M\'ethodes d'alg\`ebre abstraite en
g\'eom\'etrie alg\'ebrique, Ergebnisse der math., Springer-
Verlag, 1955.

\bibitem{T1} J.H. Teh, A homology and cohomology for real projective
varieties, Preprint in arXiv.org, math.AG/0508238.

\bibitem{T2} J.H. Teh, Harnack-Thom Theorem for higher cycle
groups, Preprint in arXiv.org, math.AG/0509149.

\bibitem{Thom} R. Thom, Sur l'homologie des vari\`et\`es
alg\`ebriques r\`eelles, Differential and combinatorial topology
(A symposium in honor of Marston Morse), University Press,
Princeton, N.J., 1965, 255-265.
\end{thebibliography}

\end{document}